\documentclass[aop,citesort,MSNbibl,dvips]{arximspdf}
\usepackage{graphicx}

\doi{10.1214/10-AOP541}
\volume{38}
\issue{6}
\pubyear{2010}
\firstpage{2418}
\lastpage{2442}

\makeatletter

\newcommand{\iint}{\int\!\!\int}

\newtheorem{theorem}{Theorem}
\newtheorem{cor}[theorem]{Corollary}
\newtheorem{prop}[theorem]{Proposition}
\newtheorem{lem}[theorem]{Lemma}
\newproclaim{stassump}{Standing assumption}

\makeatother

\begin{document}
\begin{frontmatter}

\title{Curvature, concentration and error estimates for Markov chain
Monte Carlo}
\runtitle{Curvature, concentration and MCMC error estimates}

\begin{aug}
\author[A]{\fnms{Ald\'{e}ric} \snm{Joulin}\ead[label=e1]{ajoulin@insa-toulouse.fr}} and
\author[B]{\fnms{Yann} \snm{Ollivier}\corref{}\ead[label=e2]{yann.ollivier@umpa.ens-lyon.fr}}
\runauthor{A. Joulin and Y. Ollivier}
\affiliation{Universit\'{e} de Toulouse and CNRS and \'{E}cole normale
sup\'{e}rieure de Lyon}
\address[A]{Universit\'{e} de Toulouse\\
Institut National des Sciences Appliqu\'{e}es\\
Institut de Math\'{e}matiques de Toulouse\\
31077 Toulouse\\
France\\
\printead{e1}} 
\address[B]{UMPA\\
\'{E}cole normale sup\'{e}rieure de Lyon\\
46, all\'{e}e d'Italie\\
69007 Lyon\\
France\\
\printead{e2}}
\end{aug}

\received{\smonth{4} \syear{2009}}
\revised{\smonth{2} \syear{2010}}

%
\begin{abstract}
We provide explicit nonasymptotic estimates for the rate of convergence
of empirical means of Markov chains, together with a Gaussian or
exponential control on the deviations of empirical means. These
estimates hold under a ``positive curvature'' assumption expressing a
kind of metric ergodicity, which generalizes the Ricci curvature from
differential geometry and, on finite graphs, amounts to contraction
under path coupling.
\end{abstract}

%

\end{frontmatter}

\setcounter{section}{-1}
\section{Introduction}\label{intro}
The goal of the Markov chain Monte Carlo method is to provide an efficient
way to approximate the integral $\pi(f):=\int f(x) \pi({
{d}}x)$ of a
function $f$
under a finite measure $\pi$ on some space $\mathcal{X}$. This approach,
which has been very successful, consists in constructing a hopefully
easy-to-simulate Markov chain
$(X_1,X_2,\ldots,X_k,\ldots)$ on $\mathcal{X}$ with stationary distribution
$\pi$,
waiting for a time $T_0$ (the \textit{burn-in}) so that the chain gets
close to its stationary distribution,
and then estimating $\pi(f)$ by the empirical mean on the next $T$ steps
of the trajectory, with $T$ large enough:
\[
\hat\pi(f):=\frac1T \sum_{k=T_0+1} ^{T_0+T} f(X_k).
\]
We refer, for example, to \cite{RR} for a review of the topic.

Under suitable assumptions \cite{MT}, it is known that
$\hat\pi(f)$ almost surely tends
to $\pi(f)$ as $T\to\infty$, that
the variance of
$\hat\pi(f)$ decreases asymptotically like $1/T$ and that a
central limit theorem holds for the errors $\hat\pi(f)-\pi(f)$.
Unfortunately, these theorems are asymptotic only, and thus mainly of
theoretical interest since they do not allow to
give explicit confidence intervals for $\pi(f)$ at a given time $T$. Some
even say that confidence intervals disappeared the day MCMC methods
appeared.

In this paper, we aim at establishing rigorous nonasymptotic upper
bounds for the error $| \hat\pi(f)-\pi(f) |$, which will
provide good deviation estimates and confidence intervals for $\pi(f)$.
An important\vadjust{\goodbreak} point is that we will try to express all results in terms of
explicit quantities that are readily computable given a choice of a
Markov chain; and, at the same time, recover correct order of magnitudes
and improve on known estimates in a surprising variety of examples.

Our nonasymptotic estimates have the same qualitative behavior as theory
predicts in the asymptotic regime: the variance of $\hat\pi(f)$
decreases like $1/T$, and the bias
decreases exponentially in $T_0$.
Moreover, we provide a Gaussian or exponential control on
deviations of $\hat\pi(f)$, which allows for good confidence
intervals.
Finally, we find
that the
influence of the choice of the starting point on the variance of
$\hat\pi(f)$ decreases like $1/T^2$.

Our results hold under an assumption of \textit{positive curvature}
\cite{Oll07,Oll09}, which can be understood as a kind of ``metric
ergodicity'' expressing contraction in a transportation distance. This
assumption reduces to the well-known contraction under path coupling when
$\mathcal{X}$ is a finite graph (see, e.g., Chapter 14 in \cite
{LPW09} or
\cite{BD97}), and, when
$\mathcal{X}$ is a manifold, to positive (Ricci) curvature in the
ordinary geometric
sense.

Not all ergodic Markov chains satisfy this assumption, but important
examples include spin systems at high temperature, several of the usual
types\vspace*{1pt} of waiting queues, processes such as the Ornstein--Uhlenbeck
process on ${\mathbb{R}}^d$ or Brownian motion on positively curved
manifolds. We
refer to \cite{Oll09} for more examples and discussions on how one can
check this assumption, but let us stress out that, at least in principle,
this curvature can be computed explicitly given a Markov transition
kernel. This property or similar ones using contraction in
transportation distance can be traced back to Dobrushin
\cite{Dob70,DS85}, and have appeared several times independently in the
Markov chain literature
\cite{CW94,Dob96,BD97,DGW04,Jou07,Oll07,Oll09,Jou,Oli}.

Similar concentration inequalities have been recently investigated in
\cite{Jou} for time-continuous Markov jump processes. More precisely, the
first author obtained Poisson-type tail estimates for Markov processes
with positive Wasserstein curvature. Actually, the latter
is nothing but a continuous-time version of the Ricci curvature
emphasized in the present paper, so that we expect to recover such
results by a simple limiting argument (cf. Section \ref{sect:examples}).

Gaussian-like estimates for the deviations of empirical means have
previously been given in \cite{Lez98} using the \textit{spectral gap} of
the Markov chain, under different conditions (namely that the chain is
reversible, that the law of the initial point has a density w.r.t.
$\pi$, and that $f$ is bounded). The positive curvature assumption, which
is a stronger property than the spectral gap used by Lezaud, allows to
lift these restrictions: our results apply to an arbitrary starting
point, the function $f$ only has to be Lipschitz, and reversibility plays
no particular role. In a series of papers (e.g.,
\cite{Wu00,CG08,GLWY}), the spectral approach has been extended into a
general framework
for deviations of empirical means using various types of functional
inequalities; in particular \cite{GLWY} contains a very nice
characterization of asymptotic variance of empirical means of Lipschitz
functions in terms of a functional inequality $W_1I$ satisfied by the
invariant distribution.\vadjust{\goodbreak}

\section{Preliminaries and statement of the results}\label{sec1}

\subsection{Notation}

\subsubsection*{Markov chains} In this paper, we consider a Markov chain
$(X_N)_{N\in{\mathbb{N}}}$ in a Polish (i.e., metric, complete, separable)
state space $(\mathcal{X}, d)$. The associated transition kernel is denoted
$(P_x)_{x\in\mathcal{X}}$ where each $P_x$ is a probability measure
on $\mathcal{X}$, so
that $P_x({{d}}y)$ is the transition probability from $x$ to
$y$. The
$N$-step transition kernel is defined inductively
as
\[
P_x^N({{d}}y):=
\int_\mathcal{X}P_x ^{N-1} ({{d}}z) P_z ({{d}}y)
\]
(with $P_x^1:=P_x$).
The distribution at time $N$ of the Markov chain given the initial
probability measure $\mu$ is the measure $\mu P^N$ given by
\[
\mu P ^N ({{d}}y)
=
\int_\mathcal{X}P_x ^N ({{d}}y) \mu({{d}}x).
\]
Let as usual $\mathbb{E}_x$ denote the expectation of a random
variable knowing that the initial point of the Markov chain is $x$.
For any measurable function $f\dvtx \mathcal{X}\to{\mathbb{R}}$, define
the iterated averaging operator as
\[
P^N f(x):=\mathbb{E}_x f(X_N)=\int_\mathcal{X}f(y) P_x ^N ({
{d}}y),\qquad x\in\mathcal{X}.
\]

A probability measure $\pi$ on $\mathcal{X}$ is said to be \textit
{invariant} for
the chain if $\pi= \pi P$. Under suitable assumptions on the Markov
chain $(X_N)_{N\in{\mathbb{N}}}$, such an invariant measure $\pi$
exists and is
unique, as we will see below.

Denote by $\mathcal{P}_d (\mathcal{X})$ the set of those probability
measures $\mu$ on
$\mathcal{X}$
such that $\int_\mathcal{X}d(y,x_0) \mu({{d}}y) <\infty$
for some (or
equivalently for all) $x_0 \in\mathcal{X}$. We will always assume
that the map
$x\mapsto P_x$ is measurable, and that $P_x\in\mathcal{P}_d(\mathcal
{X})$ for every
$x\in
\mathcal{X}$. These assumptions are always satisfied in practice.

\subsubsection*{Wasserstein distance}
The $L^1$ transportation distance, or
Wasserstein distance, between two probability measures $\mu_1 , \mu_2
\in\mathcal{P}_d (\mathcal{X})$ represents the ``best'' way to send
$\mu_1$ on $\mu_2$
so that on average, points are moved by the smallest possible distance.
It is defined \cite{Vil03} as
\[
W_1(\mu_1 , \mu_2) := \inf_{\xi\in\Pi(\mu_1, \mu_2)} \int
_\mathcal{X}
\int_\mathcal{X}d(x,y) \xi({{d}}x,{{d}}y) ,
\]
where $\Pi(\mu_1, \mu_2)$ is the set of probability measures $\xi$
on $\mathcal{P}
_d ( \mathcal{X}\times\mathcal{X})$ with marginals $\mu_1$ and $\mu
_2$, that is, such
that $\int_y \xi({{d}}x,{{d}}y)=\mu_1({
{d}}x)$ and $\int_x \xi({{d}}x,{{d}}
y)=\mu_2({{d}}y)$. [So, intuitively, $\xi({
{d}}x,{{d}}y)$ represents the
amount of
mass traveling from $x$ to $y$.]

\subsubsection*{Ricci curvature of a Markov chain}
Our main assumption in this paper is the following, which can be seen
geometrically as a
``positive Ricci curvature'' \cite{Oll09} property of the Markov chain.
\begin{stassump*}
There exists $\kappa>0$ such that
\[
W_1(P_x,P_y)\leq(1-\kappa) \,d(x,y)
\]
for any
$x,y\in\mathcal{X}$.
\end{stassump*}

When the space $\mathcal{X}$ is a finite graph, this is equivalent to the
well-known \textit{path coupling} criterion; for example, in Theorem
14.6 of
\cite{LPW09}, the coefficient $1-\kappa$ appears as $e^{-\alpha}$.

In practice, 
it is not
necessary to compute the exact value of the Wasserstein distance
$W_1(P_x,P_y)$: it is enough to exhibit one choice of $\xi({
{d}}x,{{d}}y)$
providing a good value of $W_1(P_x,P_y)$.

An important remark is that on a ``geodesic'' space $\mathcal{X}$,
it is sufficient to control $W_1(P_x,P_y)$ only for
nearby points $x,y \in\mathcal{X}$, and not for all pairs of points
(Proposition 19 in \cite{Oll09}). For instance, on a graph, it is enough
to check the assumption on pairs of neighbors.

These remarks make the assumption possible to check in practice, as we
will see from the examples below.

\subsubsection*{More notation: Eccentricity, diffusion constant, local
dimension, granularity} Under the assumption above, Corollary 21 in
\cite{Oll09} entails the existence of a unique invariant measure $\pi
\in\mathcal{P}_d (\mathcal{X})$ with, moreover, the following
geometric ergodicity in
$W_1$-distance [instead of the classical total variation distance,
which is obtained by choosing the trivial metric $d(x,y) = 1 _{\{ x\neq
y \}}$]:
%
%
\begin{equation}
\label{eq:T_1_cv}
W_1(\mu P^N , \pi) \leq(1-\kappa) ^N W_1(\mu, \pi) ,
\end{equation}
and in particular
%
%
\begin{equation}
\label{eq:T_1_cv_eccent}
W_1(P_x ^N , \pi) \leq(1-\kappa) ^N E(x) ,
\end{equation}
where the \textit{eccentricity} $E$ at point $x\in\mathcal{X}$ is
defined as
\[
E(x)
:= \int_\mathcal{X}d(x,y) \pi({{d}}y).
\]
Eccentricity will play the role that the diameter
of the graph plays in the path coupling
setting, though on unbounded spaces better bounds are needed. Indeed,
eccentricity satisfies the
bounds \cite{Oll09}
\[
E(x) \leq
\cases{\operatorname{diam}\mathcal{X};\cr
E(x_0)+d(x,x_0) , &\quad $x_0 \in\mathcal{X}$; \cr
\displaystyle\frac1\kappa\int_\mathcal{X}d(x,y) P_x({{d}}y).}
\]
These a priori estimates are useful in various situations. In
particular, the last one is ``local'' in the sense that it is
easily computable given the Markov kernel $P_x$.

Let us also introduce the \textit{coarse diffusion
constant} $\sigma(x)$ of the Markov chain at a point $x\in\mathcal
{X}$, which
controls the size of the steps, defined
by
\[
\sigma(x)^2:=\frac12 \iint d(y,z)^2 P_x({{d}}y) P_x({{d}}z).
\]
Let the \textit{local dimension} $n_x$ at point $x\in\mathcal{X}$ be
given by
\begin{eqnarray*}
n_x:\!&=&
\mathop{\inf_{f\dvtx \mathcal{X}\to{\mathbb{R}}}}
_{f\ 1\mbox{-}\mathrm{Lipschitz}}
\frac{
\iint d(y,z)^2 P_x({{d}}y) P_x({{d}}z)
}{\iint| f(y)-f(z) |^2 P_x({{d}}y) P_x({{d}}z)}
\\
&\geq& 1.
\end{eqnarray*}
Let the \textit{granularity} of the Markov chain be
\[
\sigma_\infty:=\frac12 \sup_{x\in\mathcal{X}} \operatorname
{diam}\operatorname{Supp}P_x,
\]
which we will often assume to be finite.

For example, for the simple random walk in a graph we have
$\sigma_\infty\leq1$ and $\sigma(x)^2\leq2$. The bound $n_x\geq1$ is
often sufficient for application to graphs.

Finally, we will denote by $\|\cdot\|_{\mathrm{Lip}}$ the
usual Lipschitz seminorm of a function $f$ on $\mathcal{X}$:
\[
\|f \|_{\mathrm{Lip}} := \sup_{x\neq y} \frac{| f(x)-f(y) |}{d(x,y)}.
\]

\subsection{Results}

Back to the \hyperref[intro]{Introduction}, choosing integers $T\geq1$ and $T_0\geq0$
and setting
\[
\hat\pi(f):=\frac1T \sum_{k=T_0+1} ^{T_0+T} f(X_k),
\]
the purpose of this paper is to understand how fast the difference
$
| \hat\pi(f) -\pi(f) |
$
goes to 0 as $T$ goes to infinity for a large class of functions $f$.
Namely, we will consider Lipschitz functions (recall that Corollary 21 in
\cite{Oll09}
implies that all
Lipschitz functions are $\pi$-integrable).

\subsubsection*{Bias and nonasymptotic variance}
Our first interest is in
the nonasymptotic mean quadratic error
\[
\mathbb{E}_x [ | \hat\pi(f) -\pi(f) | ^2 ]
\]
given any starting point $x\in\mathcal{X}$ for the
Markov chain.

%
%
\begin{figure}

\includegraphics{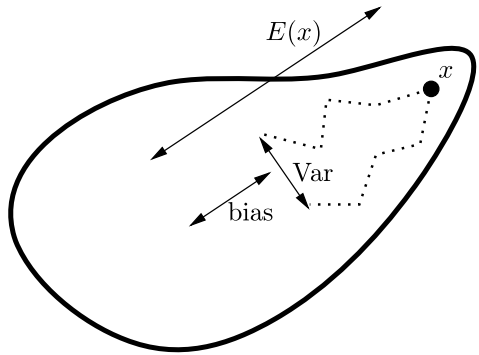}

\caption{Bias and variance.}
\label{fig1}
\end{figure}

There are two contributions to this error (see Figure \ref{fig1}): a \textit{variance} part,
controlling how $\hat\pi(f)$ differs between two independent runs
both starting at $x$, and a \textit{bias} part, which is the difference
between $\pi(f)$ and the average value of $\hat\pi(f)$
starting at $x$. Namely, the mean quadratic error decomposes as the
sum of the squared bias plus the variance
%
%
\begin{equation}
\label{eq:decompo}
\mathbb{E}_x [ | \hat\pi(f) -\pi(f) |^2 ]=
| \mathbb{E}_x \hat\pi(f)-\pi(f) |^2
+
\operatorname{Var}_x \hat\pi(f),
\end{equation}
where $\operatorname{Var}_x \hat\pi(f):=\mathbb{E}_x [| \hat\pi
(f) -\mathbb{E} _x \hat\pi(f) |^2 ]$.

As we will see, these two terms have different behaviors
depending on $T_0$ and~$T$. For instance, the bias is expected to
decrease exponentially fast as the burn-in period $T_0$ is large, whereas
if $T$ is fixed, the variance term does not vanish as $T_0 \to\infty$.

Let us start with control of the bias term, which depends, of course, on
the starting point of the Markov chain.
All proofs are postponed to Section~\ref{sect:proof}.

\begin{prop}[(Bias of empirical means)]
\label{prop:bias}
For any Lipschitz function $f\dvtx\mathcal{X}\to{\mathbb{R}}$, we have
the upper bound on
the bias
%
%
\begin{equation}
\label{eq:bias}
| \mathbb{E}_x \hat\pi(f) - \pi(f) |\leq\frac{(1-\kappa)^{T_0
+1}}{\kappa T} E(x) \|f \|_{\mathrm{Lip}}.
\end{equation}
\end{prop}

The variance term is more delicate to control. For comparison,
let us first
mention that under the invariant measure $\pi$, the variance of a
Lipschitz function $f$ is bounded as follows (and this estimate is
often sharp):
%
%
\begin{equation}
\label{eq:var_invar}
\operatorname{Var}_\pi f \leq\|f \|_{\mathrm{Lip}}^2 \sup_{x\in
\mathcal{X}} \frac{\sigma(x)^2}{ n_x
\kappa}
\end{equation}
(cf. Lemma \ref{lemme:var_invar} below or Proposition 32 in
\cite{Oll09}). This implies that, were one able to sample from the
invariant distribution $\pi$, the ordinary Monte Carlo method of
estimating $\pi(f)$ by the average over $T$ independent samples would
yield a variance bounded by
$
\frac{\|f \|_{\mathrm{Lip}}^2}{T} \sup_{x\in\mathcal{X}} \frac
{\sigma(x)^2}{n_x
\kappa}.
$
Because of correlations, this does not hold for the MCMC method.
Nevertheless, we get the following.
\begin{theorem}[(Variance of empirical means, 1)]
\label{thm:var1}
Provided the inequalities make sense, we have
%
%
\begin{equation}
\label{eq:var1}
\operatorname{Var}_x \hat\pi(f) \leq
\cases{\displaystyle
\frac{\|f \|_{\mathrm{Lip}}^2}{\kappa T} \sup_{x\in\mathcal{X}}
\frac{\sigma(x)^2}{n_x \kappa}, &\quad if $T_0=0$,\cr
\displaystyle\frac{\|f \|_{\mathrm{Lip}}^2}{\kappa T} \biggl( 1+\frac{1}{\kappa T} \biggr)
\sup_{x\in\mathcal{X}} \frac{\sigma(x)^2}{n_x \kappa}, &\quad otherwise.}
\end{equation}
\end{theorem}

The most important feature of this formula is the $1/(\kappa T)$ factor,
which means there is an additional $1/\kappa$ factor with respect to the
ordinary Monte Carlo case. Intuitively, the idea is that
correlations disappear after roughly
$1/\kappa$ steps, and so $T$ steps of the MCMC method are ``worth''
only $\kappa T$
independent samples. This $1/(\kappa T)$ factor will appear repeatedly in
our text.

To get convinced that this $1/(\kappa T)$ factor in our estimate (\ref
{eq:var1}) is natural, observe
that if the burn-in $T_0$ is large enough, then the law of $X_{T_0}$ will
be close to the invariant distribution $\pi$ so that $\operatorname
{Var}_x \hat
\pi(f)$
will behave like $\operatorname{Var}_{X_0\sim\pi} \hat\pi(f)$.
Then we have
\begin{eqnarray*}
&&\operatorname{Var}_{X_0\sim\pi} \hat\pi(f)\\
&&\qquad=\frac{1}{T^2} \Biggl(\sum_{i=1}^T
\operatorname{Var}_{X_0\sim\pi}(f(X_i))+2\sum_{1\leq i<j\leq
T}\operatorname{Cov}_{X_0\sim\pi
}(f(X_i),f(X_j)) \Biggr),
\end{eqnarray*}
but our assumption on $\kappa$ easily implies that correlations
decrease exponentially fast with rate $1-\kappa$ so that (at least in the
reversible case) we have $\operatorname{Cov}_{X_0\sim\pi
}(f(X_0),f(X_t))\leq
(1-\kappa)^t\operatorname{Var}_{X_0\sim\pi} f(X_0)$. In particular,
for any fixed $i$,
we have $\operatorname{Var}_{X_0\sim\pi}
(f(X_i))+2\sum_{j>i}\operatorname{Cov}_{X_0\sim\pi
}(f(X_i),f(X_j))\leq\break
\frac{2}{\kappa}
\operatorname{Var}_{X_0\sim\pi} f(X_i)$. Plugging this into the
above
yields $\operatorname{Var}_{X_0\sim\pi} \hat\pi(f)\leq\frac
{2}{\kappa T}\operatorname{Var}
_\pi
f$, which explains the $1/(\kappa T)$ factor.

\subsubsection*{Unbounded diffusion constant}
In the formulas above, a supremum of $\sigma(x)^2/ n_x $ appeared. This
is fine when
considering, for example, the simple random walk on a graph, because then
$\sigma(x)^2/n_x\approx1$ for all $x\in\mathcal{X}$. However, in
some situations
(e.g., binomial distributions on the cube), this supremum is
much larger than a typical value, and, in some continuous-time limits on
infinite spaces, the supremum may even be infinite (as in the example
of the $M/M/\infty$ queueing process below). For such situations, one
expects the variance to depend, asymptotically, on the average of
$\sigma(x)^2/n_x$ under the invariant
measure $\pi$, rather than its supremum.

The next result\vspace*{1pt} generalizes Theorem \ref{thm:var1} to Markov chains with
unbounded diffusion constant $\sigma(x)^2/n_x$. We will assume that
$\sigma(x)^2/n_x$ has at most linear growth (this is consistent with the
usual theorems on Markov processes, in which linear growth on the
coefficients of the diffusion is usually assumed).

Of course, if one starts the Markov chain at a point $x$ with large
$\sigma(x)^2$, meaning that the chain has a large\vadjust{\goodbreak}
diffusion constant at the
beginning, then the variance of the empirical means started
at $x$ will be accordingly large, at least for small $T$. This gives
rise, in the estimates below, to a variance term depending on
$x$; it so happens that this terms decreases
like $1/T^2$ with time.
\begin{theorem}[(Variance of empirical means, 2)]
\label{thm:var2}
Assume that there exists a Lipschitz function $S$ with $\|S \|_{\mathrm
{Lip}} \leq
C$ such that
\[
\frac{\sigma(x)^2}{n_x\kappa} \leq S(x) ,\qquad x\in\mathcal{X}.
\]
Then the variance of the empirical mean is bounded as follows:
%
%
\begin{equation}
\label{eq:var2}
\operatorname{Var}_x \hat\pi(f) \leq
\cases{\displaystyle
\frac{\|f \|_{\mathrm{Lip}}^2}{\kappa T} \biggl(\mathbb{E}_\pi S +
\frac{C}{\kappa T} E (x) \biggr), \qquad \mbox{if $T_0=0$},\vspace*{2pt}\cr
\displaystyle\frac{\|f \|_{\mathrm{Lip}}^2}{\kappa T} \biggl( \biggl( 1+\frac{1}{\kappa T}
\biggr) \mathbb{E}_\pi S + \frac{2C (1-\kappa)^{T_0} }{\kappa T} E (x)
\biggr),\vspace*{2pt}\cr
\hspace*{146.4pt}\mbox{otherwise}.}
\end{equation}
\end{theorem}

In particular, the upper bound behaves asymptotically like
$\|f \|_{\mathrm{Lip}}^2 \mathbb{E}_\pi S /(\kappa
T)$, with a correction of order $1/(\kappa T)^2$ depending on the
initial point $x$.

Note that in some situations, $\mathbb{E}_\pi S$ is known in advance from
theoretical reasons. In general, it is always possible to chose any
origin $x_0\in\mathcal{X}$ (which may or may not be the initial point
$x$) and
apply the estimate
$\mathbb{E}_\pi S \leq S(x_0 )+CE(x_0)$.

\subsubsection*{Concentration results}
To get good confidence intervals for $\hat\pi(f)$, it is
necessary to investigate deviations, that is, the behavior of the
probabilities
\[
{\mathbb{P}}_x \bigl( | \hat\pi(f) - \pi(f) | > r \bigr) ,
\]
which reduce to deviations of the centered empirical mean
\[
{\mathbb{P}}_x \bigl( | \hat\pi(f) - \mathbb{E}_x \hat\pi(f) | > r
\bigr),
\]
if the bias is known. Of course,
the Bienaym\'{e}--Chebyshev inequality states that
${\mathbb{P}}_x ( | \hat\pi(f) - \mathbb{E}_x \hat\pi(f) | > r
)\leq\frac{\operatorname{Var}_x \hat\pi(f)}{r^2}$, but this does
not decrease
very fast with~$r$, and Gaussian-type deviation estimates are
often necessary to get good confidence intervals.
Our next results
show that the probability of a deviation of size $r$ is bounded by an
explicit
Gaussian or exponential term. (The same Gaussian-exponential transition
also appears in \cite{Lez98} and other works.)

Of course, deviations for the function $10f$ are $10$ times bigger
than deviations for $f$, so we will use the rescaled
deviation $\frac{\hat\pi(f) - \mathbb{E}_x
\hat\pi(f)}{\|f \|_{\mathrm{Lip}}}$.

In the sequel, we assume that $\sigma_\infty<\infty$. Once more the
proofs of the following results are established in
Section \ref{sect:proof}.
\begin{theorem}[(Concentration of empirical means, 1)]
\label{thm:conc1}
Denote by $V^2$ the quantity
\[
V^2 := \frac{1}{\kappa T} \biggl( 1+ \frac{T_0}{T} \biggr) \sup_{x\in\mathcal
{X}} \frac
{\sigma(x) ^2}{n_x \kappa} .
\]
Then empirical means satisfy the following concentration result:
%
%
\begin{equation}
\label{eq:conc1}
{\mathbb{P}}_x \biggl( \frac{| \hat\pi(f) - \mathbb{E}_x \hat\pi(f)
|}{\|f \|_{\mathrm{Lip}}} \geq r \biggr) \leq\cases{
2 e^{-{r^2}/({16V^2})}, &\quad if $r \in(0, r_{\max} )$,\cr
2 e^{- {\kappa T r}/({12 \sigma_\infty})}, &\quad if $r \geq
r_{\max}$,}
\end{equation}
where the boundary of the Gaussian window is given by $ r_{\max} :=
4V^2 \kappa T/ (3\sigma_\infty)$.
\end{theorem}

Note that $r_{\max}\gg1/\sqrt{T}$ for large $T$, so that the
Gaussian window gets better and better when normalizing by the standard
deviation. This is in accordance with the central limit theorem
for Markov chains \cite{MT}.

As for the case of variance above, we also provide an estimate using
the average of
$\sigma(x)^2/n_x$ rather than its supremum.
\begin{theorem}[(Concentration of empirical means, 2)]
\label{thm:conc2}
Assume that there exists a Lipschitz function $S$ with $\|S \|_{\mathrm
{Lip}} \leq
C$ such that
\[
\frac{\sigma(x)^2}{n_x\kappa} \leq S(x),\qquad x\in\mathcal{X}.
\]
Denote by $V_x ^2$ the following term depending on the initial
condition $x$:
\[
V_x ^2 := \frac{1}{\kappa T} \biggl( 1+\frac{T_0}{T} \biggr) \mathbb{E}
_\pi S + \frac{C E (x)}{\kappa^2 T^2}.
\]
Then the following concentration inequality holds:
%
%
\begin{equation}
\label{eq:conc2}
{\mathbb{P}}_x \biggl( \frac{| \hat\pi(f) - \mathbb{E}_x \hat\pi(f)
|}{\|f \|_{\mathrm{Lip}}}\geq r \biggr) \leq\cases{
2 e^{- {r^2}/({16 V_x ^2})}, &\quad if $r \in(0, r_{\max} )$, \cr
2 e^{- {\kappa Tr}/({4 \max\{ 2C,3\sigma_\infty\}})}, &\quad if $r
\geq
r_{\max}$,}\hspace*{-28pt}
\end{equation}
where $r_{\max} := 4V_x ^2 \kappa T /\max\{ 2C,3\sigma_\infty\}$.
\end{theorem}

The two quantities $V^2$ and $V_x ^2$ in these theorems are essentially
similar to the estimates of the empirical
variance $\operatorname{Var}_x \hat\pi(f)$ given in Theorems \ref{thm:var1}
and \ref{thm:var2}, so that the same comments apply.

\subsubsection*{Randomizing the starting point} As can be seen from the
above,
we are mainly concerned
with Markov chains starting from a deterministic point.
The random case might be treated as follows. Assume that the starting
point $X_0$ of the chain is taken at random according to some probability
measure $\mu$. Then an\vadjust{\goodbreak} additional variance term appears in the
variance/bias decomposition (\ref{eq:decompo}), namely
\begin{eqnarray*}
\mathbb{E}_{X_0\sim\mu} [ | \hat\pi(f) -\pi(f) |^2 ] & = &
| \mathbb{E}_{X_0\sim\mu} \hat\pi(f) - \pi(f) |^2
+
\int_\mathcal{X}\operatorname{Var}_x \hat\pi(f) \mu({{d}}x)
\\
& &{} + \operatorname{Var}[ \mathbb{E}( \hat\pi(f) | X_0 ) ].
\end{eqnarray*}
The new variance term depends on how ``spread'' the initial distribution
$\mu$ is and can be easily bounded. Indeed, we have
\[
\mathbb{E}\bigl( \hat\pi(f) | X_0=x \bigr) = \frac1T \sum_{k= T_0 +1}
^{T_0 +T} P^k f(x)
\]
so that if $f$ is, say, $1$-Lipschitz,
\begin{eqnarray*}
\operatorname{Var}[ \mathbb{E}( \hat\pi(f) | X_0 ) ]
& = &
\frac{1}{2T^2 } \int_\mathcal{X}\int_\mathcal{X}\Biggl| \sum_{k= T_0
+1} ^{T_0 +T} \bigl(P^k f(x)- P^k f(y) \bigr) \Biggr|^2 \mu({{d}}x) \mu
({{d}}y)
\\
& \leq& \frac{(1-\kappa)^{2(T_0 +1)} }{ 2\kappa^2 T^2} \int
_\mathcal{X}
\int
_\mathcal{X}
d(x,y)^2 \mu({{d}}x) \mu({{d}}y) ,
\end{eqnarray*}
since $P^kf$ is $(1-\kappa)^k$-Lipschitz. This is fast-decreasing both
in $T_0$ and $T$.

Note also that the bias can be significantly reduced if $\mu$ is
known, for some reason, to be close to the invariant measure $\pi$. More
precisely, the eccentricity $E(x)$ in the bias formula (\ref{eq:bias})
above is replaced with the $L^1$ transportation distance $W_1(\mu,\pi)$.

\subsubsection*{Comparison to spectral methods} Previous results on
deviations of empirical means often rely on spectral methods (e.g.,
\cite{Lez98,Rud}). When the starting point of the MCMC method is taken
according to an initial distribution $\mu$, a factor
$\operatorname{Var}_\pi({{{d}}\mu}/{{{d}}\pi})$
generally appears
in the deviation estimates
\cite{Lez98,GLWY}, the same way it does for convergence of $P^N$ to the
invariant distribution \cite{DS96}. In particular, when the initial
distribution $\mu$ is a Dirac measure (as may be the case for practical
reasons), these estimates can perform poorly
since $\operatorname{Var}_\pi({{{d}}\delta_{x_0}}/{{
{d}}\pi})$
behaves like the cardinality of
$\mathcal{X}$ and is not even finite on continuous spaces $\mathcal
{X}$. (This was one of
the motivations for the introduction of logarithmic Sobolev inequalities
on discrete spaces, see \cite{DS96}.) This forces either to start with a
distribution $\mu$ close enough to $\pi$---but estimating $\pi$ is often
part of the problem; or to use a nonzero burn-in time $T_0$ so that the
distribution at time $T_0$ is close enough to $\pi$. But then
$T_0$ has to be comparable to or larger than the mixing time of the
chain, and since estimations of the mixing time by spectral methods have
exactly the same shortcomings, another ingredient (e.g., logarithmic
Sobolev inequalities) is needed to efficiently bound~$T_0$.
On the other hand, the approach
used here performs well when starting at a single point, even with small
or vanishing burn-in time.

\subsubsection*{Convergence of $\hat\pi$ to $\pi$} The fact that
the empirical measure $\hat\pi$ yields estimates close to
$\pi
$ when integrating
Lipschitz\vadjust{\goodbreak} functions does not mean that $\hat\pi$ itself is
close to
$\pi$. To see this, consider the simple case when $\mathcal{X}$ is a
set of $N$
elements equipped with any metric. Consider the trivial Markov chain
on $\mathcal{X}$ which sends every point $x\in\mathcal{X}$ to the
uniform probability
measure on $\mathcal{X}$ (so that $\kappa=1$ and the MCMC method
reduces to the
ordinary Monte
Carlo method). Then it is clear that for any function $f$, the bias
vanishes and the empirical variance is
\[
\operatorname{Var}_x \hat\pi(f)=\frac1T \operatorname{Var}_\pi f,
\]
which in particular does not depend directly on $N$ and allows to
estimate $\pi(f)$ with a sample size independent of $N$, as is
well known to any
statistician. But
the empirical measure $\hat\pi$ is a sum of Dirac masses at $T$
points, so that its Wasserstein distance to the uniform measure cannot
be small
unless
$T$ is comparable to $N$.

This may seem to contradict the Kantorovich--Rubinstein duality
theorem \cite{Vil03},
which states that
\[
W_1(\hat\pi,\pi)=\sup_{f\ 1\mbox{-}\mathrm{Lipschitz}}
\hat\pi(f)-\pi(f).
\]
Indeed, we know that for a function $f$ fixed in advance, very
probably $\hat\pi(f)$ is close to $\pi(f)$. But for every
realization of the random measure $\hat\pi$ there may be a
particular function $f$ yielding a large error. What is true, is that
the \textit{averaged} empirical measure $\mathbb{E}_x \hat\pi$
starting at
$x$ tends to
$\pi$ fast enough, namely
\[
W_1(\mathbb{E}_x \hat\pi,\pi)\leq
\frac{(1-\kappa)^{T_0
+1}}{\kappa T} E(x),
\]
which is just a restatement of our bias
estimate above (Proposition \ref{prop:bias}). But as we have just seen,
$\mathbb{E}_x
W_1(\hat\pi,\pi)$ is generally much larger.

\section{Examples and applications}
\label{sect:examples}

We now show how these results can be applied to various settings where
the positive curvature assumption is satisfied, ranging from discrete
product spaces to waiting queues, diffusions on ${\mathbb{R}}^d$ or
manifolds, and
spin systems. In several examples, our results improve on the
literature.

\subsection{A simple example: Discrete product spaces}
Let us first consider a very simple example. This is
mainly illustrative, as in this case the invariant measure $\pi$ is
very easy to simulate.
Let $\mathcal{X}=\{0,1\}^N$ be the space of $N$-bit sequences equipped
with the uniform
probability measure. We shall use the Hamming distance on $\mathcal
{X}$, namely,
the distance between two sequences of $0$'s and $1$'s is the number of
positions at which they differ. The Markov chain we shall consider
consists, at each step, in choosing a position $1\leq i\leq N$ at random,
and replacing the $i$th bit of the sequence with either a $0$ or a $1$
with probability $1/2$. Namely, starting at $x=(x_1,\ldots,x_N)$, we have
$P_x(x)=1/2$ and $P_x(x_1,\ldots,1-x_i,\ldots,x_N)=1/2N$.

A typical Lipschitz function for this example is the function $f_0$
equal to the proportion of ``$0$'' bits in the sequence, for which
$\|f_0 \|_{\mathrm{Lip}}=1/N$.

Then an elementary computation (Example 8 in \cite{Oll09}) shows that
$\kappa=1/N$, so that our theorems apply. The various quantities of
interest are estimated as $\sigma_\infty=1$, $\sigma(x)^2\leq2$ and
$n_x\geq1$; using Remark 40 in \cite{Oll09} yields a slightly better
estimate $\frac{\sigma(x)^2}{n_x}\leq1/2$. Moreover, $E(x)=N/2$ for any
$x\in\mathcal{X}$.\vspace*{2pt}

Then our bias estimate (\ref{eq:bias}) for a Lipschitz function $f$ is
\[
| \mathbb{E}_x \hat\pi(f) - \pi(f) |\leq
\frac{N^2}{2T} (1-1/N)^{T_0+1} \|f \|_{\mathrm{Lip}}
\leq\frac{N^2}{2T} e^{-T_0/N} \|f \|_{\mathrm{Lip}}.
\]
So taking
$T_0\approx2N\log N$ is enough to ensure small bias. This estimate of the
mixing time is known to be the correct order of magnitude:
indeed, if each bit has been updated at least once
(which occurs after a time $\approx N\log N$) then the measure is
exactly the invariant measure and so, under this event, the bias exactly
vanishes.
In contrast, the classical
estimate using only the spectral gap yields only $O(N^2)$ for the
mixing time
\cite{DS96}.

The variance estimate (\ref{eq:var1}) reads
\[
\operatorname{Var}\hat\pi(f)\leq\frac{N^2}{2T}(1+N/T) \|f \|
_{\mathrm{Lip}}^2
\]
so that, for example, for the function $f_0$ above, taking $T\approx N$
will yield a variance \mbox{$\approx$}$1/N$, the same order of magnitude as the
variance of $f_0$ under the uniform measure. (With a little work, one can
convince oneself that this order of magnitude is correct for large $T$.)

The concentration result (\ref{eq:conc1}) reads, say with $T_0=0$, and for
the Gaussian part
\[
{\mathbb{P}}_x \bigl( | \hat\pi(f) - \mathbb{E}_x \hat\pi(f) | \geq r
\bigr) \leq2
e^{-Tr^2/8N^2\|f \|_{\mathrm{Lip}}^2}
\]
so that, for example, for $f_0$ we simply get get $2e^{-Tr^2/8}$. For comparison,
starting at a Dirac and without burn-in
the spectral estimate from \cite{Lez98} behaves like
$2^{N/2}e^{-Tr^2/4}$ for
small $r$,
so that we roughly improve the estimate by a factor $2^{N/2}$ in the
absence of burn-in, due to the
fact that the density of the law of
the starting point (a Dirac mass) plays no role in our setting, as
discussed above in the comparison with spectral methods.

\subsection{Heat bath for the Ising model}
Let $G$ be a finite graph. Consider the classical Ising model from
statistical mechanics \cite{Mar04}, namely the
configuration space $\mathcal{X}:=\{-1,1\}^G$ together with the energy function
$U(s):=-\sum_{x\sim y \in G} s(x)s(y)-h\sum_{x\in G} s(x)$ for $s\in
\mathcal{X}
$, where
$h\in{\mathbb{R}}$. For some $\beta\geq0$, equip
$\mathcal{X}$
with the Gibbs distribution $\pi:=e^{-\beta U} /Z$ where as usual
$Z:=\sum_s e^{-\beta U(s)}$. The distance between two states is defined
as the number of vertices of $G$ at which their values differ, namely
$d(s,s'):=\frac12\sum_{x\in G}| s(x)-s'(x) |$.

For $s\in\mathcal{X}$ and $x\in G$, denote by $s_{x+}$ and $s_{x-}$
the states
obtained from $s$ by setting $s_{x+}(x)=+1$ and $s_{x-}(x)=-1$,
respectively. Consider the following random walk on $\mathcal{X}$,
known as the
\textit{heat bath} or \textit{Glauber dynamics} \cite{Mar04}: at each step,
a vertex $x\in G$ is chosen at
random, and a new value for $s(x)$ is picked according to local
equilibrium, that is, $s(x)$ is set to $1$ or $-1$ with probabilities
proportional to $e^{-\beta U(s_{x+})}$ and $e^{-\beta U(s_{x-})}$,
respectively (note that only the neighbors of $x$ influence the ratio of
these probabilities). The Gibbs distribution $\pi$ is invariant (and
reversible) for this
Markov chain.

When $\beta=0$, this Markov chain is identical to the Markov chain on
$\{0,1\}^N$ described above, with $N=| G |$. Therefore, it comes as no
surprise that for $\beta$ small enough, curvature is positive. More
precisely, one finds \cite{Oll09}
\[
\kappa\geq\frac{1}{| G |} \biggl(1-v_{\max} \frac{e^\beta
-e^{-\beta
}}{e^{\beta}+e^{-\beta}} \biggr),
\]
where $v_{\max}$ is\vspace*{1pt} the maximal valency of a vertex of $G$. In
particular, if
$
\beta< \frac12 \ln
(\frac{v_{\max}+1}{v_{\max}-1} )
$
then $\kappa$ is positive.
This is not surprising, as the current research interest in
transportation distances can be traced back to \cite{Dob70} (where the
name \textit{Vasershtein distance} is introduced), in which a
criterion for convergence of spin systems is introduced.
Dobrushin's criterion was a contraction property of the Markov chain in
Wasserstein distance, and thus, in this context, precisely coincides with
our notion of $\kappa>0$. (See also \cite{Per}.)

Let us see\vspace*{1pt} how our theorems apply, for example, to the magnetization
$f_0(s):=\frac{1}{| G |}\sum_{x\in G} s(x)$. With the metric we use,
we have
$\|f_0 \|_{\mathrm{Lip}}=\frac{2}{| G |}$.

Let
$\gamma:=1-v_{\max} \frac{e^\beta-e^{-\beta}}{e^{\beta
}+e^{-\beta}}$,
so that $\kappa=\frac{\gamma}{| G |}$,
and assume that $\gamma>0$.
Using the gross inequalities $\sigma_\infty=1$, $\sigma(s)^2\leq2$ and
$n_s\geq1$, the
variance estimate of Theorem~\ref{thm:var1} reads, with $T_0=0$,
\[
\operatorname{Var}_s \hat\pi(f_0) \leq\frac{8}{\gamma^2 T},
\]
where $s$ is any initial configuration.
For example, taking $T\approx| G |$ (i.e., each site of $G$ is
updated a few times by the heat bath) ensures that $\operatorname{Var}_s
\hat\pi(f_0)$ is of the same order of magnitude as the variance of
$f_0$ under the invariant measure.

Theorem \ref{thm:conc1} provides a Gaussian estimate for deviations, with
similar variance up to numerical constants.
The transition for Gaussian to non-Gaussian regime becomes
very relevant when the external magnetic field $h$ is large, because then
the number of spins opposing the magnetic field has a Poisson-like rather
than Gaussian-like behavior (compare Section 3.3.3 in \cite{Oll09}).

The bias is controlled as follows: using $E(s)\leq\operatorname{diam}
\mathcal{X}=| G |$ in Proposition \ref{prop:bias} one finds
$| \mathbb{E}_s \hat\pi(f_0)-\pi(f_0) |\leq
2| G |(1-\gamma/| G |)^{T_0}/\gamma T$ so that taking
${T_0\approx| G |\log}| G |$ is a good choice.

These results are not easily compared with the literature, which often
focusses on getting nonexplicit constants for systems of infinite size
\cite{Mar04}. However, we have seen that even in the case $\beta=0$ our
estimates improve on the spectral estimate, and our results provide very
explicit bounds on the time necessary to run a heat bath simulation, at
least for $\beta$ not too large.

\subsection{The $M/M/\infty$ queueing process}
We now focus on a continuous-time example, namely the $M/M/\infty$
queueing process. This is a continuous-time Markov chain $(X_t)_{t\geq
0}$ on ${\mathbb{N}}$ with transition kernel given for small $t$ by
\[
P_t (x,y) = \cases{
\lambda t + o(t), &\quad if $y=x+1$, \cr
x t + o(t), &\quad if $y=x-1$, \cr
1- (\lambda+x )t + o(t), &\quad if $y=x$,}
\]
where $\lambda$ is a positive parameter.
The (reversible) invariant
measure is the Poisson distribution $\pi$ on ${\mathbb{N}}$ with parameter
$\lambda$. Although this process is very simple in appearance, the
unboundedness of the associated transition rates makes the determination
of concentration inequalities technically challenging. Here, we will
get a
convenient concentration inequality for Lipschitz functions $f$ with
respect to the classical metric on ${\mathbb{N}}$, in contrast with
the situation
of \cite{Jou} where Poisson-like concentration estimates are provided for
Lipschitz functions with respect to an ad hoc metric. The techniques
used here allow us to overcome this difficulty.

First, let us consider, given $d\in{\mathbb{N}^\star}$, $d >\lambda
$, the
so-called binomial Markov chain $(X_N ^{(d)})_{N\in{\mathbb{N}}}$ on
$\{ 0,1,
\ldots, d\}$, with transition probabilities given by
\[
P_x ^{(d)} (y) = \cases{
\displaystyle\frac{\lambda}{d} \biggl( 1- \frac{x}{d} \biggr), &\quad if $y=x+1$;\cr
\displaystyle\biggl( 1- \frac{\lambda}{d} \biggr) \frac{x}{d}, &\quad if $y=x-1$; \cr
\displaystyle\frac{\lambda x}{d^2} + \biggl( 1- \frac{\lambda}{d} \biggr) \biggl( 1- \frac{x}{d}
\biggr), &\quad if
$y=x$.}
\]
The invariant measure is the binomial distribution $\pi^{(d)}$ on $\{
0,1, \ldots, d\}$ with parameters $d$ and $\lambda/d$. It is not
difficult to show that the Ricci curvature is $\kappa=1/d$ and that
$\sigma(x) ^2\leq(\lambda+x)/d$ for $x\in\{ 0,1, \ldots, d\}$.

But now,
take instead the continuous-time version of the above, namely the
Markov process $(X_t ^{(d)})_{t\geq0}$ whose transition kernel is
defined for any $t\geq0$ as
\[
P_t ^{(d)}(x,y) = e^{-t} \sum_{k =0} ^{+\infty} \frac{t^k}{k!} \bigl(P_x
^{(d)} \bigr)^k (y),\qquad x, y \in\{ 0,1, \ldots, d\} .
\]
As $d\to\infty$, the invariant measure $\pi^{(d)}$ converges weakly to
the Poisson measure~$\pi$, which is nothing but the invariant measure of
the $M/M/\infty$ queueing process. One can check (using, e.g., Theorem 4.8
in \cite{JS}) that the process $(X_t ^{(d)})_{t\geq0}$ sped up by a
factor $d$ converges to the $M/M/\infty$ queueing process
$(X_t)_{t\geq
0}$ in a
suitable sense (in the \textit{Skorokhod space} of
c\`adl\`ag functions equipped with the Skorokhod topology).

To derive a concentration inequality for the empirical mean
$\hat\pi(f) := t^{-1} \times\break\int_0 ^t f(X_s) \,{{d}}s$, where $f$ is
$1$-Lipschitz on ${\mathbb{N}}$ and time $t$ is fixed, we proceed as follows.
First, we will obtain a concentration estimate for the continuous-time
binomial Markov chain $(X_t ^{(d)})_{t\geq0}$ by using Theorem \ref
{thm:conc2} for the chain
$(X_{\varepsilon N} ^{(d)})_{N\in{\mathbb{N}}}$ with $\varepsilon
\to0$, and then we will
approximate the $M/M/\infty$
queueing process $(X_t)_{t\geq0}$ by the sped-up process $(X_{td}
^{(d)})_{t\geq0}$ with $d\to\infty$.

For small $\varepsilon$, the Markov chain $(X_{\varepsilon N}
^{(d)})_{N\in{\mathbb{N}}}$ has
Ricci curvature bounded below by $\varepsilon/d$, eccentricity $E(x)
\leq
x+E(0)=x+\lambda$, square diffusion constant $\sigma(x) ^2$ of
order $\varepsilon(\lambda+x)/d$, and $n_x\geq1$, so that we may
take $S(x)
:= \lambda+x$ in
Theorems \ref{thm:var2} and \ref{thm:conc2} above (with $T_0 = 0$
for simplicity). Let $f$ be a $1$-Lipschitz function.
For a given $t> 0$, we have ${\mathbb{P}}
_x$-almost surely the Riemann approximation
\[
\hat\pi^{(d)} (f) := \frac{1}{t} \int_0 ^t f \bigl(X_s ^{(d)}\bigr) \,{
{d}}s =
\lim_{T \to+\infty} \hat\pi^{(d),T} (f),
\]
where $\hat\pi^{(d),T} (f) := \frac{1}{T} \sum_{k=1} ^T f
(X_{kt/T} ^{(d)} )$. So applying Theorem \ref{thm:conc2} to the
Markov chain $(X_{\varepsilon N} ^{(d)})_{N\in{\mathbb{N}}}$ with
$\varepsilon=t/T$, we get by
Fatou's lemma
\begin{eqnarray*}
&&
{\mathbb{P}}_x \bigl( \bigl| \hat\pi^{(d)} (f) - \mathbb{E}_x \hat\pi^{(d)}
(f) \bigr| >r \bigr)\\
&&\qquad \leq
\liminf_{T \to+\infty} {\mathbb{P}}_x \bigl( \bigl| \hat\pi^{(d),T} (f) -
\mathbb{E}_x \hat\pi^{(d),T} (f) \bigr| >r \bigr)
\\
&&\qquad \leq
\cases{
2 e^{- {t^2 r^2}/({16 d (2\lambda t+ (\lambda+x)d)})}, &\quad if $r \in
\bigl(0, r_{\max} ^{(d)}\bigr)$,\cr
2 e^{- {tr}/({12d})}, &\quad if $r \geq r_{\max} ^{(d)}$,}
\end{eqnarray*}
where $r_{\max} ^{(d)} := ( 8\lambda t+ 4(\lambda+x)d
)/3t$. Finally, we approximate $(X_t)_{t\geq0}$ by the sped-up
process $(X_{td} ^{(d)})_{t\geq0}$ with $d\to\infty$ and apply Fatou's
lemma again to obtain
the following.
\begin{cor}
Let $(X_s)_{s\geq0}$ be the $M/M/\infty$ queueing process with
parameter $\lambda$.
Let $f\dvtx{\mathbb{N}}\to{\mathbb{R}}$ be a $1$-Lipschitz function.
Then for any $t>0$, the
empirical mean $\hat\pi(f):=t^{-1}\int_{s=0}^t f(X_s)\, {
{d}}s$ under
the process starting at
$x\in{\mathbb{N}}$ satisfies the
concentration inequality
\[
{\mathbb{P}}_x \bigl( | \hat\pi(f) - \mathbb{E}_x \hat\pi(f) | >r \bigr)
\leq
\cases{2 e^{- {t r^2}/({16 (2\lambda+
(\lambda+x)/t)})}, &\quad if $r \in(0, r_{\max} )$,\cr
2 e^{- {tr}/{12}}, &\quad if $r \geq r_{\max}$,}
\]
where $r_{\max} := ( 8\lambda t+ 4(\lambda+x) )/3t$.
\end{cor}

Let us mention that a somewhat similar, albeit much less explicit,
concentration inequality has been derived in \cite{GLWY} via
transportation-information inequalities and a drift condition of
Lyapunov-type.

Our results generalize to other kinds of waiting queues, such as queues
with a finite number of servers and positive abandon rate.

\subsection{Euler scheme for diffusions} Let
$(X_t)_{t\geq0}$ be the solution of the following stochastic
differential equation on the Euclidean space ${\mathbb{R}}^d$:
\[
{{d}}X_t = b(X_t) \,{{d}}t + \sqrt{2} \rho(X_t)\,
{{d}}W_t,
\]
where $(W_t)_{t\geq0}$ is a standard Brownian motion in ${\mathbb
{R}}^d$, the
function $b\dvtx{\mathbb{R}}^d \to{\mathbb{R}}^d$ is measurable, as is
the $d\times d$
matrix-valued function $\rho$.
For a given matrix $A$, we define the Hilbert--Schmidt
norm $\|A \|_{\mathrm{HS}} := \sqrt{\operatorname{tr} AA^*}$ and the operator
norm $\|A \|_{{\mathbb{R}}^d}:=\sup_{v\neq0}\frac{\|Av \|}{\|v \|}$.

We assume that the following stability condition \cite{BHW97,DGW04} is
satisfied:
\begin{itemize}
\item[(C)]
the functions $b$ and $\rho$ are Lipschitz, and
there exists $\alpha>0$ such that
\[
\|\rho(x) - \rho(y) \|^2_{\mathrm{HS}}+\langle x-y , b(x)-b(y)
\rangle\leq- \alpha\|x-y \|^2 ,\qquad x,y \in{\mathbb{R}}^d.
\]
\end{itemize}
A typical example is the Ornstein--Uhlenbeck process, defined by
$\rho=\operatorname{Id}$ and $b(x)=-x$. As we will see, this
assumption implies that
$\kappa>0$.

The application of Theorems \ref{thm:var2}, \ref{thm:conc1} and
\ref{thm:conc2} on this example requires careful approximation arguments
(see below). The result is the following.
\begin{cor}\label{corrA}
Let $\hat\pi(f):=t^{-1}\int_{s=0}^t f(X_s)\, {{d}}s$
be the empirical
mean of the $1$-Lipschitz function $f\dvtx{\mathbb{R}}^d\to{\mathbb{R}}$
under the diffusion process
$(X_t)_{t\geq0}$ above, starting at point $x$. Let
$S\dvtx\mathcal{X}\to{\mathbb{R}}$ be a $C$-Lipschitz function with
$S(x)\geq\frac
{2}{\alpha}
\|\rho(x) \|_{{\mathbb{R}}^d}^2$. Set
\[
V^2_x:=\frac{1}{\alpha t} \mathbb{E}_\pi S+ \frac{CE(x)}{\alpha^2 t^2}.
\]
Then
one has
$\operatorname{Var}_x \hat\pi(f)\leq V^2_x$ and
\[
{\mathbb{P}}_x \bigl( | \hat\pi(f) - \mathbb{E}_x \hat\pi(f) |
> r \bigr) \leq\cases{
2 e^{- {r^2}/({16 V_x ^2})}, &\quad if $r \in
(0, r_{\max} )$,\cr
2 e^{- {\alpha tr}/({8C})}, &\quad if $r \geq r_{\max}$,}
\]
where $r_{\max} := 2V_x ^2 \alpha t /C$.
\end{cor}

An interesting\vspace*{-1pt} case is when $\rho$ is constant or bounded, in which case
one can take $S(x):=\sup_x \frac{2 \|\rho(x) \|_{{\mathbb
{R}}^d}^2}{\alpha}$ so
that $C=0$. Then $r_{\max}=\infty$
and the exponential regime disappears. For this particular case, our
result is comparable to \cite{GLWY}, but note, however, that their result
requires some regularity on the distribution of the starting point of the
process, in contrast with ours.

Note that the final result features the average of $S$ under the
invariant distribution. Sometimes this value is known from theoretical reasons,
but in any case the assumption $(C)$ implies very explicit bounds on the
expectation of $d(x_0,x)^2$ under the invariant measure $\pi$ \cite{BHW97},
which can be used to bound $\mathbb{E}_\pi S$ knowing $S(x_0)$, as
well as to
bound $E(x)$.

The Lipschitz growth of $\|\rho\|_{{\mathbb{R}}^d}^2$ allows to
treat stochastic
differential equations where the diffusion constant $\rho$ grows like
$\sqrt{x}$, such as naturally appear in population dynamics or
superprocesses.
\begin{pf*}{Proof of Corollary \ref{corrA}}
Consider the underlying Euler scheme with (small) constant step
$\delta t$ for the stochastic differential equation above, that is, the
Markov chain $(X_{N} ^{(\delta t)})_{N\in{\mathbb{N}}}$
defined by
\[
X_{N+1} ^{(\delta t)} = X_{N} ^{(\delta t)} + b\bigl(X_{N} ^{(\delta t)}\bigr)
\delta t+ \sqrt{2\delta t} \rho\bigl(X_{N} ^{(\delta t)}\bigr) Y_N,
\]
where $(Y_N)$ is any sequence of i.i.d. standard Gaussian random
vectors. When $\delta t\to0$, this process tends to a weak solution of
the stochastic differential equation~\cite{BHW97}.

Let us see how Theorems \ref{thm:var2}, \ref{thm:conc1} and
\ref{thm:conc2} may be applied. The measure $P_x$ is
a Gaussian
with expectation $x+b(x)\delta t$ and covariance matrix $2\delta t\rho
\rho^* (x)$.
Let $(X_{N} ^{(\delta t)}(x))_{N\in{\mathbb{N}}}$ be the chain
starting at $x$.
Under (C), we have
\begin{eqnarray*}
\mathbb{E}\bigl[ \bigl\Vert X_{1} ^{(\delta t)}(x)-X_{1} ^{(\delta t)}(y)\bigr\Vert
^2 \bigr] & =
& \Vert x-y\Vert^2 + 2\delta t\langle x-y , b(x)-b(y) \rangle
\\& &{}
+ 2 \delta t\Vert\rho(x) - \rho(y) \Vert_{\mathrm{HS}} ^2
+ \delta t^2 \Vert b(x)-b(y) \Vert^2 \\
& \leq& \Vert x-y\Vert^2 \bigl( 1 - \alpha\delta t+ O(\delta t^2) \bigr) ^2,
\end{eqnarray*}
so that we obtain $\kappa\geq\alpha\delta t+ O(\delta t^2)$.
Moreover,
the diffusion constant $\sigma(x) $ is given by
\begin{eqnarray*}
\sigma(x)^2 & = &\frac{1}{2} \int_{{\mathbb{R}}^d} \int_{{\mathbb
{R}}^d} \|y-z \|^2
P_x ({{d}}y) P_x ({{d}}z) \\
& = & 2\delta t\Vert\rho(x) \Vert_{\mathrm{HS}} ^2
\end{eqnarray*}
by a direct computation.

Next, using the Poincar\'{e} inequality for Gaussian measures in
${\mathbb{R}}^d$,
with a little work one gets that the local dimension is
\[
n_x= \frac{\|\rho(x) \|_{\mathrm{HS}} ^2}{\|\rho(x) \|_{{\mathbb{R}}^d}^2}.
\]
For
example, if $\rho$ is the $d\times d$ identity matrix we have $n_x=d$,
whereas $n_x=1$ if $\rho$ is of rank $1$.

So, we get that $\frac{\sigma(x)^2}{\kappa n_x}$ is bounded by the
function
\[
S(x) := \frac{2}{\alpha} \|\rho(x) \|_{{\mathbb{R}}^d}^2
+O(\delta t).
\]

However, here we have $\sigma_\infty=\infty$. This can be circumvented
either by directly plugging into Lemma \ref{lem:laplace} the
well-known Laplace
transform estimate for Lipschitz functions of Gaussian variables, or
slightly changing the approximation scheme as
follows. Let us assume that $\sup_x \|\rho(x) \|_{{\mathbb
{R}}^d}<\infty$. Now,
replace the Gaussian random vectors $Y_N$ with random vectors whose law
is supported in a large ball of radius $R$ and approximates a Gaussian
(the convergence theorems of \cite{BHW97} cover this situation as well).
Then we have $\sigma_\infty=R\sqrt{2\delta t} \sup_x \|\rho(x) \|
_{{\mathbb{R}}^d}$.
This modifies the quantities $\sigma(x)^2$ and $\rho(x)$ by a factor at
most $1+o(1)$ as $R\to\infty$.

Therefore, provided $S$ is Lipschitz,
we can apply Theorem \ref{thm:conc2} to the empirical mean
$\hat\pi(f):=t^{-1}\int_{s=0}^t f(X_s) \,{{d}}s$ by using the Euler
scheme at time $T=t/\delta t$ with $\delta t\to0$, and using Fatou's lemma
as we did above for the case of the $M/M/\infty$ process.
Note in particular that $\sigma_\infty\to0$ as $\delta t\to0$, so that
$\sigma_\infty$ will disappear from $r_{\max}$ in the final result.

Finally, the constraint $\sup_x \|\rho(x) \|_{{\mathbb{R}}^d}<\infty
$ can be lifted
by considering that, under our Lipschitz growth assumptions on $b$ and
$\|\rho(x) \|_{{\mathbb{R}}^d}^2$, with arbitrary high probability
the process does
not leave a compact set and so, up to an arbitrarily small error, the
deviation probabilities considered depend only on the behavior of $\rho$
and $b$ in a compact set.
\end{pf*}

\subsection{Diffusions on positively curved manifolds}
Consider a diffusion process $(X_t)_{t\geq0}$ on a smooth, compact
$N$-dimensional Riemannian manifold $M$, given by the
stochastic differential equation
\[
{{d}}X_t=b \,{{d}}t+\sqrt{2} \,{{d}}B_t
\]
with infinitesimal generator
\[
L:=\Delta+b\cdot\nabla,
\]
where $b$ is a vector field on $M$, $\Delta$ is the Laplace--Beltrami
operator and $B_t$ is the standard Brownian motion in the Riemannian
manifold $M$. The Ricci curvature
of this operator in the Bakry--\'{E}mery sense \cite{BE85}, applied to
a tangent vector
$v$, is
$\operatorname{Ric}(v,v)-v\cdot\nabla_{ v} b$ where $\operatorname
{Ric}$ is the usual Ricci tensor.
Assume that this quantity is at least
$K$ for any unit tangent vector $v$.

Consider as above the Euler approximation scheme at time $\delta t$ for this
stochastic differential equation: starting at a point $x$, follow the
flow of $b$ for
a time $\delta t$, to obtain a point $x'$; now take a random tangent
vector $w$ at $x'$ whose law is a Gaussian in the tangent plane at $x'$ with
covariance matrix equal to the metric, and follow the geodesic generated
by $w$ for a time $\sqrt{2\delta t}$. Define $P_x$ to be the law of the
point so obtained. When $\delta t\to0$, this Markov chain approximates
the process $(X_t)_{t\geq0}$ (see, e.g., Section I.4 in \cite{Bis81}).
Just as above, actually the
Gaussian law has to be truncated to a large ball so that
$\sigma_\infty<\infty$.

For this Euler approximation, we have $\kappa\geq K\delta t
+O(\delta t^{3/2})$
where $K$ is a lower bound for Ricci--Bakry--\'{E}mery
curvature \cite{Oll09}. We have $\sigma(x)^2=2N\delta t+O(\delta t
^{3/2})$ and
$n_x=N+O(\sqrt{\delta t})$. The details are omitted, as they are very
similar to the case of ${\mathbb{R}}^d$ above, except that in a
neighborhood of
size $\sqrt{\delta t}$ of a given point, distances are distorted by a
factor $1\pm O(\sqrt{\delta t})$ w.r.t. the Euclidean case. We restrict
the statement to compact manifolds so that the constants hidden in the
$O(\cdot)$ notation are uniform in $x$.

So applying Theorem \ref{thm:conc1} to the Euler scheme at time
$T=t/\delta t$ we get the following corollary.
\begin{cor}
Let $(X_t)_{t\geq0}$ be a process as above on a smooth, compact
$N$-dimensional Riemannian manifold
$\mathcal{X}$,
with Bakry--\'{E}mery curvature at least $K>0$.
Let $\hat\pi(f):=t^{-1}\int_{s=0}^t f(X_s) \,{{d}}s$ be the empirical
mean of the $1$-Lipschitz function $f\dvtx\mathcal{X}\to{\mathbb{R}}$
under the diffusion process
$(X_t)$ starting at some point $x\in\mathcal{X}$. Then
\[
{\mathbb{P}}_x\bigl(| \hat\pi(f)-\mathbb{E}_x \hat\pi(f) | > r\bigr)\leq
2e^{-{K^2
t r^2}/{32}}.
\]
\end{cor}

Once more, a related estimate appears in \cite{GLWY}, except that
their result features an additional factor $\|{{d}}\beta
/{{d}}\pi\|_2$ where
$\beta$ is the law of the initial point of the Markov chain and $\pi$
is the invariant
distribution, thus preventing it from being applied with $\beta$ a Dirac
measure at $x$.

\subsection{Nonlinear state space models}
Given a Polish state space $(\mathcal{X}, d)$, we consider the Markov chain
$(X_N)_{N\in{\mathbb{N}}}$ solution of the following equation:
\[
X_{N+1} = F(X_N , W_{N+1}),\qquad X_0 \in\mathcal{X},
\]
which models a noisy dynamical system. Here, $(W_N)_{N\in{\mathbb
{N}}}$ is a
sequence of i.i.d. random variables with values in some parameter
space, with common distribution~$\mu$.
We assume that
there exists some $r<1$ such that
%
%
\begin{equation}
\label{eq:djguwu}
\mathbb{E}d (F(x,W_1) ,F(y,W_1) ) \leq r d(x,y) ,\qquad x,y \in\mathcal{X},
\end{equation}
and that moreover the following function is $L^2$-Lipschitz on
$\mathcal{X}$:
\[
x\mapsto\mathbb{E}[ d(F(x,W_1) ,F(x,W_2) ) ^2 ] .
\]
Note that the assumption (\ref{eq:djguwu}) already appears in
\cite{DGW04} [condition (3.3)]
to study the propagation of Gaussian concentration to path-dependent
functionals of $(X_N)_{N\in{\mathbb{N}}}$.

Since the transition probability $P_x$ is the image measure of $\mu$ by
the function $F(x, \cdot)$, it is straightforward that the Ricci
curvature $\kappa$ is at least $1-r$, which is positive. Hence, we may
apply Theorem \ref{thm:var2} with the $L^2/(2(1-r))$-Lipschitz function
\[
S(x) := \frac{1}{2(1-r)} \mathbb{E}[ d(F(x,W_1) ,F(x,W_2) ) ^2 ] ,
\]
to obtain the variance inequality
\[
\sup_{\|f \|_{\mathrm{Lip}} \leq1} \operatorname{Var}_x \hat\pi
(f) \leq
\cases{
\displaystyle\frac{1}{(1-r) T} \biggl\{ \frac{L^2}{2(1-r)^2 T} E (x) + \mathbb{E}_\pi S
\biggr\}, & \quad if
$T_0 = 0$,
\vspace*{2pt}\cr
\displaystyle\frac{1}{\kappa T} \biggl\{ \biggl( 1+\frac{1}{(1-r)T} \biggr) \mathbb{E}_\pi S\vspace*{2pt}\cr
\hspace*{33.7pt}{} +
\displaystyle\frac{L^2 r
^{T_0} }{(1-r)^2 T} E (x) \biggr\},
&\quad otherwise.}
\]
Note that to obtain a qualitative concentration estimate via
Theorem \ref{thm:conc2}, we need the additional assumption $\sigma
_\infty<\infty$, which depends on the properties of $\mu$ and of the
function $F(x, \cdot)$ and states that at each step the noise has a bounded
influence.

\section{Proofs}
\label{sect:proof}

\subsection{\texorpdfstring{Proof of Proposition
\protect\ref{prop:bias}}{Proof of Proposition 1}}

Let $f$ be a $1$-Lipschitz function. Let us recall from \cite{Oll09} that
for $k\in{\mathbb{N}}$, 
the function $P^k f$ is $(1-\kappa)^k$-Lipschitz.
Then we have by the invariance of $\pi$
\begin{eqnarray*}
| \mathbb{E}_x \hat\pi(f) - \pi(f) | & = &
\frac1T \Biggl| \sum_{k= T_0 +1} ^{T_0 +T} \int_\mathcal{X}\bigl( P^k f(x) -
P^k f(y) \bigr) \pi({{d}}y) \Biggr| \\
& \leq& \frac1T \sum_{k= T_0 +1} ^{T_0 +T} (1-\kappa) ^k \int
_\mathcal{X}
d(x,y) \pi({{d}}y) \\
& \leq& \frac{(1-\kappa) ^{T_0 +1}}{\kappa T} \int_\mathcal
{X}d(x,y) \pi
({{d}}
y) ,
\end{eqnarray*}
so that we obtain the result.

\subsection{\texorpdfstring{Proof of Theorems \protect\ref{thm:var1} and
\protect\ref{thm:var2}}{Proof of Theorems 2 and 3}}
\label{sect:proof_var}

Let us start with a variance-type result under the measure after $N$
steps. The proof relies on a simple induction argument and is left to
the reader.
\begin{lem}
\label{lemme:var_invar}
For any $N\in{\mathbb{N}^\star}$ and any Lipschitz function $f$ on
$\mathcal{X}$, we have
%
%
\begin{equation}
\label{eq:var_Mn}
P^N (f^2)- (P^N f)^2 \leq\|f \|_{\mathrm{Lip}} ^2 \sum_{k=0} ^{N-1}
(1-\kappa)
^{2(N-1-k)} P^k \biggl( \frac{\sigma^2}{n} \biggr) .
\end{equation}
In particular, if the rate $x\mapsto\sigma(x) ^2 / n_x$ is bounded,
then letting $N$ tend to infinity above entails a variance estimate
under the invariant measure $\pi$:
%
%
\begin{equation}
\label{eq:var_invar2}
\operatorname{Var}_\pi f \leq\|f \|_{\mathrm{Lip}} ^2 \sup_{x\in
\mathcal{X}} \frac{\sigma(x)
^2}{n_x \kappa} .
\end{equation}
\end{lem}

Now, we are able to prove the variance bounds of Theorems \ref{thm:var1}
and \ref{thm:var2}. Given a $1$-Lipschitz function $f$, consider the functional
\[
f_{x_1,\ldots,x_{T-1}} (x_T) := \frac1T \sum_{k=1} ^T f(x_k) ,
\]
the others coordinates $x_1,\ldots,x_{T-1}$ being fixed. The function
$f_{x_1,\ldots,x_{T-1}}$ is $1/T$-Lipschitz, hence $1/(\kappa T)$-Lipschitz
since $\kappa\leq1$.
Moreover, for each $k\in\{ T-1, T-2, \ldots, 2 \} $, the conditional
expectation of $\hat\pi(f) $ knowing $X_1 = x_1, \ldots, X_k =
x_k$ can be written in terms of a downward induction
\[
f_{x_1,\ldots, x_{k-1}}(x_k) := \int_\mathcal{X}f_{x_1,\ldots, x_k}(x_{k+1})
P_{x_k} ({{d}}x_{k+1})
\]
and
\[
f_{\varnothing} (x_1) := \int_\mathcal{X}f_{x_1}(x_2) P_{x_1}
({{d}}x_2).
\]
By Lemma 3.2 (step 1) in \cite{Jou}, we know that $f_{x_1, \ldots,
x_{k-1}}$ is Lipschitz with constant~$s_k$, where
\[
s_k := \frac1T \sum_{j=0}^{T-k} (1-\kappa)^j \leq\frac{1}{\kappa T }.
\]
Hence, we can use the variance bound (\ref{eq:var_Mn}) with $N=1$ for
the function $f_{x_1, \ldots, x_{k-1}}$, successively for $k= T , T-1,
\ldots, 2$, to obtain
\begin{eqnarray*}
\mathbb{E}_x [ \hat\pi(f) ^2 ] & = & \int_{\mathcal{X}^T} f_{x_1,
\ldots,
x_{T-1}} (x_T) ^2 P_{x_{T-1}} ({{d}}x_T) \cdots P_{x_1}
({{d}}x_2) P_x
^{T_0 +1} ({{d}}x_1) \\
& \leq& \int_{\mathcal{X}^{T-1}} f_{x_1, \ldots, x_{T-2}} (x_{T-1}) ^2
P_{x_{T-2}} ({{d}}x_{T-1}) \cdots P_{x_1} ({{d}}x_2)
P_x ^{T_0 +1} ({{d}}
x_1) \\
& &{} + \frac{1}{\kappa^2 T^2 } P^{T_0 +T -1} \biggl(\frac{\sigma
^2}{n} \biggr) (x) \\
& \leq& \int_{\mathcal{X}^{T-2}} f_{x_1, \ldots, x_{T-3}} (x_{T-2}) ^2
P_{x_{T-3}} ({{d}}x_{T-2}) \cdots P_{x_1} ({{d}}x_2)
P_x ^{T_0 +1} ({{d}}
x_1) \\
& &{} + \frac{1}{\kappa^2 T^2 } \biggl( P^{T_0 +T -2}
\biggl(\frac{\sigma^2}{n} \biggr) (x) + P^{T_0 +T -1}
\biggl(\frac{\sigma^2}{n} \biggr) (x) \biggr) \\
& \leq& \cdots\\
& \leq& \int_\mathcal{X}f_{\varnothing} (x_1) ^2 P_x ^{T_0 +1}
({{d}}x_1) +
\frac{1}{\kappa^2 T^2 } \sum_{k=T_0 +1} ^{T_0 + T-1} P^{k} \biggl(
\frac{\sigma^2}{n} \biggr) (x) \\
& \leq& (\mathbb{E}_x \hat\pi(f) ) ^2 + \frac{1}{\kappa^2
T^2 } \sum_{k=0} ^{T_0} (1-\kappa)^{2(T_0 -k)} P^{k} \biggl(
\frac{\sigma^2}{n} \biggr) (x) \\
& &{} + \frac{1}{\kappa^2 T^2 } \sum_{k=T_0 +1} ^{T_0 + T-1} P^{k}
\biggl( \frac{\sigma^2}{n} \biggr) (x),
\end{eqnarray*}
where in the last step we applied the variance inequality (\ref
{eq:var_Mn}) to the Lipschitz function $f_{\varnothing}$, with $N=T_0 +1$.
Therefore, we get
\begin{eqnarray*}
\operatorname{Var}_x \hat\pi(f)  &\leq& \frac{1}{\kappa^2 T^2 } \Biggl(
\sum_{k=0} ^{T_0} (1-\kappa)^{2(T_0 -k)} P^{k} \biggl( \frac{\sigma
^2}{n} \biggr) (x)\\
&&\hspace*{65.1pt}{} + \sum_{k=T_0 +1} ^{T_0 + T-1} P^{k} \biggl(
\frac{\sigma^2}{n} \biggr) (x) \Biggr) .
\end{eqnarray*}
Theorem \ref{thm:var1} is a straightforward consequence of the latter
inequality. To establish Theorem \ref{thm:var2}, for instance,
(\ref{eq:var2}) in the case $T_0 \neq0$, we rewrite the above as
\begin{eqnarray*}
\operatorname{Var}_x \hat\pi(f) & \leq& \frac{1}{\kappa T^2 } \Biggl\{
\sum
_{k=0} ^{T_0} (1-\kappa)^{2(T_0 -k)} P^{k} S(x) + \sum_{k=T_0 +1}
^{T_0 + T-1} P^{k} S (x) \Biggr\} \\
& \leq& \frac{1}{\kappa T^2 } \Biggl\{ \sum_{k=0} ^{T_0} (1-\kappa
)^{2(T_0 -k)} \bigl( C W_1(P_x ^k ,\pi) + \mathbb{E}_\pi S \bigr) \\
&&\hspace*{59.7pt}{} + \sum_{k=T_0 +1} ^{T_0 + T-1} \bigl( C W_1(P_x ^k ,\pi) + \mathbb
{E}_\pi S
\bigr) \Biggr\} \\
& \leq& \frac{1}{\kappa T^2 } \Biggl\{ \biggl( 1+ \frac{1}{\kappa T
} \biggr) T \mathbb{E}_\pi S + \sum_{k=0} ^{T_0} C (1-\kappa)^{2T_0 -k}
E(x) \\
&&\hspace*{106.2pt}{} + \sum_{k=T_0 +1} ^{T_0 + T-1} C(1-\kappa)^k E(x) \Biggr\} \\
& \leq& \frac{1}{\kappa T^2 } \biggl\{ \biggl( 1+ \frac{1}{\kappa T
} \biggr) T \mathbb{E}_\pi S + \frac{2C(1-\kappa)^{T_0}}{\kappa} E(x) \biggr\}.
\end{eqnarray*}
Finally, the proof in the case $T_0 =0$ is very similar and is omitted.

\subsection{\texorpdfstring{Proofs of Theorems \protect\ref{thm:conc1}
and \protect\ref{thm:conc2}}{Proofs of Theorems 4 and 5}}
\label{sect:proof_conc}

The proofs of the concentration Theorems~\ref{thm:conc1} and
\ref{thm:conc2} follows the same lines as that for variance above, except
that Laplace transform estimates $\mathbb{E}e^{\lambda f-\lambda
\mathbb{E}f}$ now play
the role of the variance $\mathbb{E}[f^2]-(\mathbb{E}f)^2$.

Assume that there exists a Lipschitz function $S\dvtx\mathcal{X}\to
{\mathbb{R}}$ with
$\|S \|_{\mathrm{Lip}} \leq
C$ such that
\[
\frac{\sigma(x)^2}{n_x\kappa} \leq S(x),\qquad x\in\mathcal{X}.
\]

Let us give first a result on the Laplace transform of Lipschitz funtions
under the measure at time $N$.
\begin{lem}
\label{lem:laplace}
Let $\lambda\in( 0, \frac{\kappa T }{\max\{ 4C,6\sigma_\infty\}}
)$. Then for any $N\in{\mathbb{N}^\star}$ and any $\frac{2}{\kappa
T}$-Lipschitz
function $f$ on $\mathcal{X}$, we have
%
%
\begin{equation}
\label{eq:laplace}
P^N (e^{\lambda f}) \leq\exp\Biggl\{ \lambda P^N f + \frac{4\lambda^2
}{\kappa T^2 } \sum_{k=0} ^{N-1} P^k S \Biggr\}.
\end{equation}
In the case $C=0$, the same formula holds for any $\lambda\in(0,
\frac
{\kappa T}{6\sigma_\infty})$.
\end{lem}
\begin{pf}
Let $f$ be $\frac{2}{\kappa T}$-Lipschitz. By Lemma 38 in \cite{Oll09},
we know that if $g$ is an $\alpha$-Lipschitz function with $\alpha
\leq
1$ and if $\lambda\in(0, \frac{1}{3\sigma_\infty})$ then we have
the estimate
\[
P ( e^{\lambda g} ) \leq\exp\{ \lambda P g + \kappa\lambda
^2\alpha
^2 S \} ,
\]
and by rescaling, the same holds for the function $f$ with $\alpha=
\frac{2}{\kappa T}$ whenever $\lambda\in(0,\frac{\kappa
T}{6\sigma_\infty})$. Moreover, the function $P^N f + \frac{4
\lambda
}{\kappa T^2 } \sum_{k=0} ^{N-1} P^k S$ is also $\frac{2}{\kappa
T}$-Lipschitz for any $N\in{\mathbb{N}^\star}$, since $\lambda\in
(0, \frac
{\kappa
T}{4C} )$. Hence, the result follows by a simple induction argument.
\end{pf}

Now let us prove Theorem \ref{thm:conc2}, using again the notation of
Section \ref{sect:proof_var} above. Theorem \ref{thm:conc1} easily follows
from Theorem \ref{thm:conc2} by taking $S := \sup_{x\in\mathcal{X}}
\frac
{\sigma(x)^2}{n_x\kappa} $ and letting $C\to0$ in the formula (\ref
{eq:conc2}).

Let $f$ be a $1$-Lipschitz function on $\mathcal{X}$ and let $\lambda
\in(
0, \frac{\kappa T }{\max\{ 4C,6\sigma_\infty\}} )$. Using the
Laplace transform estimate (\ref{eq:laplace}) with $N=1$ for the
$\frac
{2}{\kappa T}$-Lipschitz functions
\[
x_k \mapsto f_{x_1, \ldots, x_{k-1}} (x_k) + \frac{4\lambda}{\kappa
T^2} \sum_{l=0} ^{T-k-1} P^l S (x_k),
\]
successively for $k= T-1,T-2, \ldots, 2$, we have
\begin{eqnarray*}
&&\mathbb{E}_x e^{\lambda\hat\pi(f)} \\
&&\qquad = \int_{\mathcal{X}^T} e^{\lambda f_{x_1, \ldots, x_{T-1}} (x_T)}
P_{x_{T-1}} ({{d}}x_T) \cdots P_{x_1} ({{d}}x_2) P_x
^{T_0 +1} ({{d}}x_1) \\
&&\qquad \leq\int_{\mathcal{X}^{T-1}} e^{ \lambda f_{x_1, \ldots, x_{T-2}}
(x_{T-1}) + {4 \lambda^2}/({\kappa T^2 }) S(x_{T-1} ) }\\
&&\qquad\hspace*{36pt}{}\times P_{x_{T-2}}
({{d}}x_{T-1}) \cdots P_{x_1} ({{d}}x_2) P_x ^{T_0 +1}
({{d}}x_1) \\
&&\qquad \leq\int_{\mathcal{X}^{T-2}} e^{ \lambda f_{x_1, \ldots, x_{T-3}}
(x_{T-2}) + {4 \lambda^2}/({\kappa T^2 }) \sum_{l=0} ^{1} P^l
S(x_{T-2} )}\\
&&\qquad\hspace*{36pt}{}\times  P_{x_{T-3}} ({{d}}x_{T-2}) \cdots P_{x_1}
({{d}}x_2) P_x ^{T_0
+1} ({{d}}x_1) \\
&&\qquad \leq\cdots\\
&&\qquad \leq\int_\mathcal{X}e^{ \lambda f_{\varnothing} (x_1) +
{4 \lambda
^2}/({\kappa T^2 }) \sum_{l=0} ^{T-2} P^l S(x_1)} P_x ^{T_0 +1} ({{d}}x_1)
\\
&&\qquad \leq e^{\lambda\mathbb{E}_x \hat\pi(f) + {4\lambda
^2}/({\kappa T^2}) ( \sum_{l=0} ^{T-2} P^{T_0 +1 +l} S (x) + \sum_{l=0}
^{T_0} P^l S (x) )},
\end{eqnarray*}
where in the last line we applied the Laplace transform estimate (\ref
{eq:laplace}) to the $\frac{2}{\kappa T}$-Lipschitz function
\[
x_1 \mapsto f_{\varnothing} (x_1) + \frac{4 \lambda}{\kappa T^2 }
\sum
_{l=0} ^{T-2} P^l S (x_1),
\]
with $N=T_0 +1$. Therefore, we get
\[
\mathbb{E}_x e^{\lambda( \hat\pi(f) - \mathbb{E}_x \hat\pi(f) ) }
\leq
e ^{4\lambda^2 V_x ^2 } .
\]
Finally, using Markov's inequality and optimizing in $\lambda\in
( 0, \frac{\kappa T }{\max\{ 4C,6\sigma_\infty\}} )$
entails the result. This ends the proof.

\section*{Acknowledgments}
The authors would like to thank the
anonymous referee for various helpful suggestions.

%

%
\printaddresses

\end{document}